\documentclass[a4paper,12pt,reqno]{amsart}
\usepackage{amssymb}\theoremstyle{plain}

\usepackage{amscd}
\usepackage{enumitem}

\newtheorem{theorem}{Theorem}
\newtheorem*{theorem*}{Theorem}
\newtheorem{corollary}{Corollary}
\newtheorem*{corollary*}{Corollary}
\newtheorem{lemma}{Lemma}
\newtheorem*{lemma*}{Lemma}
\newtheorem{proposition}{Proposition}
\newtheorem*{proposition*}{Proposition}

\newtheorem*{conjecture*}{Conjecture}
\theoremstyle{definition}
\newtheorem{definition}{Definition}
\newtheorem*{definition*}{Definition}
\theoremstyle{remark}
\newtheorem{remark}{Remark}
\newtheorem*{remark*}{Remark}
\newtheorem{example}{Example}
\newtheorem*{problem*}{Problem}

\begin{document}

\begin{center}
\title[Ex func]{Extremal functions on moduli spaces and applications}
\end{center}

\maketitle

\begin{center}
{\bf Nikolaj M. Glazunov } \end{center}

\begin{center}
{\rm Institute of Mathematics and Informatics Bulgarian Academy of Sciences }\\
{\rm Glushkov Institute of Cybernetics NASU, Kiev, } \\

{\rm  Email:} {\it glanm@yahoo.com }
\end{center} 

\bigskip



\bigskip

\footnote{The author was supported by Simons grant 992227.}  

\begin{abstract}
Our object of study is extremal functions which are defined by distance functions of convex bodies.
 These functions take values in the moduli spaces of algebraic and geometric objects associated with these ${\mathbb Z}$-modules (geometric lattices) and with convex bodies.
In most cases, convex bodies are $2$-dimensional Minkowski balls whose  boundaries are Minkowski curves and we study lattice points on these curves.
We define and investigate extremal functions that yield  the homogeneous arithmetic minimum of a
 function in a lattice, the Hermite constant, the critical determinant 
of a body, optimal packings 
of bodies,  best  values of covering constants, and optimal 
solutions of Diophantine approximation problems.
Moreover, for two-dimensional unit Minkowski balls and Minkowski domains
 we determine  the minimal areas  of inscribed and circumscribed hexagons.
\end{abstract}


\section{Introduction}
Extremal problems and functions arise in the study of many mathematical problems: 
properties of the circle and the sphere and isoperimetric inequalities; 
Chebyshev polynomials as polynomials of least deviation from zero; homogeneous 
arithmetic minimum of a distance  function in a lattice; finding the Hermite
 constant; Minkowski's hypothesis on the critical
 determinant of a domain; Bernstein's inequality; Fourier analysis; the large 
sieve; and in others.\\
Distance functions (see section \ref{nd}) can themselves be extremal functions.\\
For example, in the well-known isoperimetric problem about a circle,
 which is still Archimedes began to study, the function $(x^2 + y^2)^\frac{1}{2}$ itself 
is an  extremal function, bounding the convex region $x^2 + y^2 < 1$ 
of the area $A = \pi$ of the curve of length $L = 2\pi$, which turns
 the isoperimetric inequality $4\pi A \le L^2$ into equality.
For a sphere, the isoperimetric problem is more complicated, but the distance 
function$(x^2 + y^2 + z^2)^\frac{1}{2}$ also plays a key role.\\
In this communication we study extremal functions, which are defined by 
distance functions  of convex bodies.  \\
One of the algebraic foundations of the Geometry of Numbers are ${\mathbb Z}$-modules (geometric lattices), which are associated with the geometric bodies studied in the problems.\\
 For a more complete study of these objects, the metric properties of the  space 
  ${\mathbb R}^n$
   are used (Minkowski's convex body theorem, and others).
  In Section \ref{mtanf} we briefly discuss the connection of our research 
with Minkowski's theory for algebraic number fields. \\
The existence of extremal functions within the framework of the study is based on the fact that in formulas (\ref{cd}), (\ref{ham}), (\ref{ha}) the infimum and supremum are achieved.\\
We construct such functions as real functions that take values in the 
moduli spaces of geometric objects associated with these convex bodies.
Examples of such functions are the homogeneous arithmetic minimum of a
 function in a lattice, the Hermite constant, and the critical determinant 
of a body.
We define and investigate extremal functions that give optimal solutions to Diophantine approximation problems,
 optimal packings of Minkowski balls and Minkowski domains, 
 and the best current 
 values of covering constants.
Moreover, for two-dimensional unit Minkowski balls and Minkowski domains
 we determine  the minimal areas  of inscribed and circumscribed hexagons.
 In Section we give applications to algebraic number theory.

 
\section{Some notations and definitions}
\label{nd}
Let $D $ be a set and $\Lambda $ be a  lattice with base 
\[
\{{\mathbf a}_1, \ldots ,{\mathbf a}_n \}
\] 
in ${\mathbb R}^n.$ \\

A lattice $\Lambda $
is {\it admissible} for body $D $ 
($ D$-{\it admissible})
if ${D} \bigcap \Lambda = \emptyset $ or $0.$\\

Let
\[
  d(\Lambda) = |\det ({\mathbf a}_1, \ldots ,{\mathbf a}_n )|
 \] 
   be the determinant of 
$\Lambda.$ \\

The infimum
$\Delta(D) $ ($ inf \; \Delta(D)$) of determinants of all lattices admissible for
$D $ is called {\em the critical determinant} 
of $D: $
\begin{equation}
\label{cd}
\Delta(D) = \inf   d(\Lambda).
\end{equation}

If there is no $D-$admissible lattices then puts
$\Delta(D) = \infty. $ \\

 A lattice 
$\Lambda $ is {\em critical}
if $ d(\Lambda) = \Delta(D).$ \\

In the real $n$-dimensional space ${\mathbb R}^n$ a distance function $F({\mathbf x})$ of variable
 ${\mathbf x} = (x_1, \ldots , x_n)$ is any function which is 
 
 \begin{enumerate}

\item non-negative, 

\item continuous, 

\item  homogeneous.

\end{enumerate} 
(1) means that $F({\mathbf x}) \ge 0$, \\
(3) means that $F(t{\mathbf x}) = tF({\mathbf x}), \; (t \ge 0, t \in {\mathbb R})$.\\

We investigate convex symmetric distance functions:\\

(4)  $F({\mathbf x}+{\mathbf y}) \le F({\mathbf x}) + F({\mathbf y})$,\\
(5) $F(-{\mathbf x}) = F({\mathbf x})$.\\


{\it Homogeneous arithmetic minimum} of the distance function in the lattice

\begin{equation}
\label{ham}
 m( F,\Lambda) = \inf   F({\mathbf a}),   {\mathbf a} \in \Lambda, \; {\mathbf a} \ne 0. 
\end{equation}
\\

{\it The Hermite's constant} of the function $F$.

\begin{equation}
\label{ha}
{\gamma}_n = \sup \frac{m( F,\Lambda)}{\{d(\Lambda)\}^\frac{1}{n}}
\end{equation}
where the supremum is taken over all lattices $\Lambda$.  \\


In the cases we are considering for a natural number $n$, such distance functions have the form

\begin{equation}
\label{df}
 F({\mathbf x})  =  ( |x_1|^p + \cdots + |x_n|^p  )^\frac{1}{p}
\end{equation}

 To the convex symmetric distance function $F({\mathbf x})$ corresponds fixed bounded symmetric about origin  convex body
  ({\it centrally symmetric convex body)}
 \begin{equation}
 \label{sb}
   F({\mathbf x}) < 1.
 \end{equation}
 
For real $p > 1$ equation (\ref{sb}) defines {\it strictly convex bodies} that do not contain {\it straight line segments on the boundary}.

 We investigate {\it strictly convex centrally symmetric bodies},
which correspond to the parameter $p > 1$, as well as the limiting cases $p=1$ and $p=\infty$
(for two-dimensional bodies of the form (\ref{sb}), the cases $p = 1$ and $p=\infty$ were already investigated by Minkowski).

  Let $x_1, \cdots x_n  \in {\mathbb R}$ and let 
 \begin{equation}
\label{smt}   D_p: \;  |x_1|^p + \cdots + |x_n|^p < 1, \; p \ge 1
\end{equation}
be $n$-dimensional (open) Minkowski's balls with boundary 
\begin{equation}
\label{mt}
{{\mathbb S}^{n-1}_c}: \;   |x_1|^c + \cdots + |x_n|^c  = 1, \;  c \ge 1.
\end{equation}
 (Minkowski's spheres).\\
 
 Below we consider the case of  $2$-dimensional (open) Minkowski's balls:
\begin{equation}
\label{s2m} 
  D_p: \;  |x|^p + |y|^p < 1, \; p > 1.
\end{equation}
 
 
    Let $D$ be a fixed bounded symmetric about origin  convex body ({\it centrally symmetric convex body} for short) with volume
$V(D)$.

\begin{proposition} 
\label{p1}
 If $D$ is symmetric about the origin and convex, then $2D$ is convex and symmetric
 about the origin.
 \end{proposition}
 
  \begin{corollary}
 \label{cor1}
 Let $m$ be integer $m \ge 0$ and $n$ be natural greater $m$. 
If $2^m D$ centrally symmetric convex body then $2^n D$ is again centrally symmetrc convex body.
\end{corollary}

\subsection{Domains.}

We consider the following classes of balls (see Section \ref{Minkowski balls}) and domains.

\begin{itemize}
  
\item  {\it  Watson domains}:    $2^m D_p$,  integer $m \ge 1$, for $1 \le p<2$;

\item {\it Davis domains}: $2^m D_p$,  integer $m \ge 1$, for $p_{0} > p \ge 2$;
 
 \item {\it Chebyshev-Mordell domains}: $2^m D_p$,  integer $m \ge 1$,  for $ p \ge p_{0}$;

\end{itemize}

\section{Minkowski's curves and algebraic number fields}
\label{mtanf}
\begin{definition}
We will call the one-dimensional Minkowski sphere (\ref{mt}) the Minkowski curve.\\
 For natural $p, d, p = 2d$ the Minkowski curves 
\begin{equation}
\label{mac} 
  {\mathbb S}^{1}_{2d} = C_{2d}: x^{2d} + y^{2d} = 1
\end{equation}
  are algebraic curves. We will call these curves Minkowski algebraic curves (algebraic curves, if it is clear from the context which curves are being considered).
\end{definition}
\begin{remark}
Over complex numbers ${\mathbb C}$ the projectivization of curve $ C_{2d}$ with  polynomial equation $ p(x,y) = x^{2d} + y^{2d} - 1 \in {\mathbb C}[x,y]$      defines the projective curve of the form $ x^{2d} + y^{2d} - z^{2d} = 0$ in the projective space
    ${\mathbb C}{\mathbb P}^2$. Its genus $g = 
    (2d - 1)(d - 1)$.
  \end{remark}
 
 \begin{proposition}
 From theorem 
  \cite{GGM:PM} (see also \cite{gl1}) we have next expressions for critical determinants and their lattices:  
\begin{enumerate} 

\item  ${\Delta^{(0)}_p} = \Delta(p, {\sigma_p}) =  \frac{1}{2}{\sigma}_{p},$ 

\item $ {\sigma}_{p} = (2^p - 1)^{1/p},$

\item  ${\Delta^{(1)}_p}  = \Delta(p,1) = 4^{-\frac{1}{p}}\frac{1 +\tau_p }{1 - \tau_p}$,   

\item  $2(1 - \tau_p)^p = 1 + \tau_p^p,  \;  0 \le \tau_p < 1.$  

\end{enumerate}
\end{proposition}
 
  For their critical lattices respectively  $\Lambda_{p}^{(0)},\; \Lambda_{p}^{(1)}$ next conditions satisfy:   $\Lambda_{p}^{(0)}$ and 
 $\Lambda_{p}^{(1)}$  are  two $D_p$-admissible lattices each of which contains
three pairs of points on the boundary of $D_p$  with the
property that 
\begin{itemize}

\item $(1,0) \in \Lambda_{p}^{(0)},$

\item $(-2^{-1/p},2^{-1/p}) \in \Lambda_{p}^{(1)},$

\end{itemize}
 (under these conditions the lattices are
uniquely defined).

\begin{example}
 Lattice $\Lambda_{p}^{(0)}$ are two-dimensional lattices in ${\mathbb R}^2$ spanned by the vectors
 \begin{itemize}
 \item[]  $\lambda^{(1)} = (1, 0),$
   \item[] $\lambda^{(2)} = (\frac{1}{2}, \frac{1}{2} \sigma_p).$
   \end{itemize}
   For lattices $\Lambda_{p}^{(1)}$, due to the cumbersomeness of the formulas for the coordinates of basis vectors of the
    lattices, as an example, we present only the lattice $\Lambda_{2}^{(1)}$.
    The lattice $\Lambda_{2}^{(1)}$ is a two-dimensional lattice in ${\mathbb R}^2$ spanned by the vectors
 \begin{itemize}
 \item[]  $\lambda^{(1)} = (-2^{-1/2},2^{-1/2}),$
   \item[] $\lambda^{(2)} = (\frac{\sqrt 6 - \sqrt 2}{4}, \frac{\sqrt 6 + \sqrt 2}{4}).$
   \end{itemize}
   \end{example}
 
\subsection{Shells of critical lattices of Minkowski curves.}

\begin{lemma}
Let $(P_x, P_y)$ be a point of the critical lattice $\Lambda$ on Minkowski curve ${\mathbb S}^{1}_{p} = C_{p}$.
Then the point $(u, v)$ that satisfies conditions
\begin{equation}
\label{es}
\left\{
                   \begin{array}{lc}
   |P_x v - P_y u| = \; d(\Lambda),\\
    |u|^p +  |v|^p = \; 1.\\
 \end{array}
                       \right.
\end{equation}
belongs to $\Lambda$ and   lies on $C_{p}$. 
The shell of points of the critical lattice $\Lambda_{p}$ on Minkowski curve $C_{p}$ contains 6 points
Point coordinates $(u, v)$ can be calculated in closed form or with any precision.
\end{lemma}
\begin{example}
Each shell of points of the critical lattices $\Lambda_{p}^{(0)}$ on Minkowski curves $C_{p}$ contain 6 points:
\[
  \pm (1, 0), 
  \]
  \[
  \left(\pm \frac{1}{2},  \pm \frac{1}{2} \sigma_p \right).
\]
  
Respectively the shell of points of the critical lattice $\Lambda_{2}^{(1)}$ on Minkowski curves $C_{p}$ contains 6 points:
\[
  (-2^{-1/2},2^{-1/2}), (2^{-1/2},-2^{-1/2}),  
  \]
  \[
  \pm \left(\frac{\sqrt 6 + \sqrt 2}{4}, \frac{\sqrt 6 - \sqrt 2}{4}\right), 
  \]
  \[
    \pm \left(\frac{\sqrt 6 - \sqrt 2}{4}, \frac{\sqrt 6 + \sqrt 2}{4} \right).
\]
\end{example}

\begin{remark}
The only integer quadratic forms of the lattices $\Lambda_{p}^{(0)}$ and $\Lambda_{p}^{(1)}$ are the quadratic forms corresponding to the lattices
$\Lambda_{2}^{(0)}$    and $\Lambda_{2}^{(1)}$.
These forms are the  form  
\begin{equation}
 x^2 - xy +y^2
 \end{equation}
 for the  lattice $\Lambda_{2}^{(0)}$ with the base $\lambda^{(1)} = (1, 0), \lambda^{(2)} = (-\frac{1}{2}, \frac{1}{2} \sqrt 3)$  and the form 
 \begin{equation}
 \label{ql21}
 x^2 -xy +y^2
\end{equation}
 for the  lattice $\Lambda_{2}^{(1)}$ with the base
  $\lambda^{(1)} = (-2^{-1/2},2^{-1/2}), \lambda^{(2)} = (\frac{\sqrt 6 + \sqrt 2}{4}, \frac{\sqrt 6 - \sqrt 2}{4})$ .
 \end{remark}

For these critical lattices, it is easy to  calculate from the corresponding theta function, their 
shells for the natural numbers.
 Corresponding theta series
   for $\Lambda_{2}^{(1)}$ with
$q=\exp(\pi i \tau)$  is defined as
\begin{equation}
 \Theta_{\Lambda_{2}^{(1)}}(\tau) = \sum_{\lambda \in \Lambda_{2}^{(1)}} q^{||\lambda||_2} 
= \sum_{x, y = -\infty}^{\infty} q^{(x^2 - xy + y^2)}  = \sum_{m=0}^{\infty} N_{m} q^{m}.
\end{equation}
 where $N_m$ is the number of times (\ref{ql21}) represents $m, N_0 = 1$.
 Thus
 \[
   \Theta_{\Lambda_{2}^{(1)}}(\tau) = 1 + 6q + 6q^{3} + 6q^{4}
    + 12q^{7} +  6q^{9} + 6q^{12} +  \cdots
     \]
 
\subsection{Strictly convex curves and smooth algebraic curves from   moduli spaces}

For $1 < p < \infty$ Minkowski curves are strictly convex curves.
\begin{remark}
The length $\ell$ of the Minkowski curve for $p > 1$ is
 \[
 {\ell} = 4 \int_0^1 (1 - y^p)^{\frac{1}{p}} dy
   \]
 \end{remark}
 
 \begin{remark}
As $\ell > 3$ by results of the Jarnic \cite{jar} the curves $C_p$ and for integer $N > 0$ their $N$-fold magnifications 
$N\cdot C_p$ about the origin for corresponding length $\ell $ contains at most 
\[
 3(4\pi)^{-\frac{1}{3}} {\ell}^{\frac{2}{3}}  + {\mathcal O}({\ell}^{\frac{1}{3}})
\]
integral lattice points.
\end{remark}
\begin{remark}
The curve $C_p$ for real $p > 1$ is three times continuity differentiable. From results of  Swinnerton–Dyer \cite{sdy} follow that the number of rational points $(\frac{m}{N},\frac{n}{N})$ on $C_p$ (the number of integral points on $N\cdot C_p$) as $N \to \infty$ is estimatrd by
\[
 |N \cdot C_p \cap {\mathbb Z}^2| \le c(C_p, \epsilon)N^{\frac{3}{5} + \epsilon}.
\]
Let $y = f(x)$ be the arc $A$ of $C_{p}$ and $f(x)$ be a transcendental analytic function on $\left[0, 1\right]$.
Then   
\[
|tA \cap {\mathbb Z}^2| \le c(f,\epsilon) t^{\epsilon}.
\]
It follows from  Bombieri-Pila formula\cite{bop}.
 \end{remark}
\subsubsection{Minkowski algebraici  curves}.
  \begin{remark}
 For algebraic curve   $C_{2d}({\mathbb C})$ its   Riemann surface is oriented sphire with $g$ handles. 
 The function field ${\mathbb F}$ on the  Riemann surface is the field of rational functions of the genus $g$. 
If two such fields are isomorphic then the corresponding algebraic curves are birationaly isomorphic.
   \end{remark}
   Let $K$ be the field of algebraic numbers and  $C_{2d}(K)$ be the projective curve over $K$. Let $\mathcal A$ be an analytic subset of $C_{2 d}({\mathbb C})$.
   By results of Bombieri-Pila \cite{bop} in the form of Gasbarri \cite{gas} we have:
  \begin{remark} 
    Let $(C_{2d}, {\mathcal L})$ be a projective variety over $K$ with ample line bundle ${\mathcal L}$ over $C_{2d}$ and let 
    $h_{\mathcal L}(\cdot)$ be a height function assosiated to it.
    Let $\mathcal A$ be a non compact Riemann surface. For relatively compact open set $U \subset {\mathcal A}$ and for
    holomorphic map $\varphi: \mathcal A \to  C_{2 d}({\mathbb C})$ with Zariski dense  set $\varphi(\mathcal A)$ and the set 
    $S_U(T) = \{ z \in U | \varphi(z) \in  C_{2 d}(K)\}$ and $h_{\mathcal L}(\varphi(z)) \le T$ for cardinality $\#S_U(T)$of the set 
    $S_U(T)$ and for any positive number $\epsilon$ we have 
    \[
    \#S_U(T) \ll \exp (\epsilon T)
    \]
    where the involved constants depend on $U, \varphi$,  and ${\mathcal L}$ but not on $T$ .
   \end{remark}

\section{Optimization in problems of Diophantine approximations, packings and coverings}

Let
 
 \begin{equation}
 \label{mdi}
|\alpha x + \beta y|^p + |\gamma x + \delta y|^p \leq c \cdot 
 |\det(\alpha \delta - \beta \gamma)|^{p/2}, 
 \end{equation} 

be a diophantine inequality defined for a given real $ p >1 $;
hear $\alpha, \beta, \gamma, \delta$ are real numbers with 
$ \alpha \delta - \beta \gamma \neq 0 .$ 

The optimization problem with respect to the 
constant $c$ in inequality (\ref{mdi}) was formulated by Minkowski \cite{Mi:DA}.

Let us formulate the minimization problem in terms of the introduced concepts.

The critical lattices of Problem (\ref{mdi}) in terms of (\ref{s2m}) are among the admissible lattices that have three pairs of points on the boundary.

Let $\{{\mathbf a}_1, {\mathbf a}_2 \}$ be the basis of an admissible lattice having 3 pairs of points on the boundary of the domain $D_p$.

Six points 
${\mathbf a}_1, {\mathbf a}_2, {\mathbf a}_1 + {\mathbf a}_2 , -{\mathbf a}_1, - {\mathbf a}_2, - ({\mathbf a}_1 +{\mathbf a}_2 )$ belong to the lattice and lie on the boundary
$D_p$.

The rough form of the optimization problem for a conditional minimum has the following representation:\\

{\it {Find the minimum}} $|det \:\{{\mathbf a}_1, {\mathbf a}_2 \}|$\\

{\it {subject to}}: ${\mathbf a}_1, {\mathbf a}_2, {\mathbf a}_1 + {\mathbf a}_2$ belong to the boundary of $D_p$.

  
\section{Moduli spaces of problems of Diophantine approximations, packings and coverings}
\label{Minkowski balls}

Below we consider centrally symmetric convex hexagons.

\subsection{Moduli spaces of two-dimensional Minkowski balls and domains}

They have the next forms:

 The Minkowski-Cohn moduli space of admissible lattices of Minkowski balls has the form
 
 \begin{equation}
 \label{mcms}
 \Delta(p,\sigma) = (\tau + \sigma)(1 + \tau^{p})^{-\frac{1}{p}}
  (1 + \sigma^p)^{-\frac{1}{p}}, 
 \end{equation}
in the domain
\begin{equation}
\label{msd}
  {\mathcal M}: \; \infty > p > 1, \; 1 \leq \sigma \leq \sigma_{p} =
 (2^p - 1)^{\frac{1}{p}},
  \end{equation}
of the $ \{p,\sigma\} $-plane, where $\sigma$ is some real parameter 

Hence from Proposition  \ref{p1} and from \cite{GGM:PM} we have
  \begin{proposition}
The Minkowski-Cohn moduli space for  admissible lattices of 
 doubled Minkowski balls $2 D_p$ (Minkowski domains $2 D_p, \; m = 1$)   has the form
 \begin{equation}
 \Delta(p,\sigma)_{2 D_p} = 4 (\tau + \sigma)(1 + \tau^{p})^{-\frac{1}{p}}
  (1 + \sigma^p)^{-\frac{1}{p}}, 
 \end{equation}
in the same domain (\ref{msd}).
 \end{proposition}
 
 If we put $m \in  {\mathbb N_0}, \; 2^0 D_p = D_p, \;   \Delta(p,\sigma)_ {D_p} = \Delta(p,\sigma)$ we can give the formula for all  moduli spaces of  Minkowski balls and domains.
 \begin{equation}
 \label{mbd}
 \Delta(p,\sigma)_{2^m D_p} = 4^m (\tau + \sigma)(1 + \tau^{p})^{-\frac{1}{p}}
  (1 + \sigma^p)^{-\frac{1}{p}}, 
 \end{equation}
 
  \begin{proposition}
\label{mpl}
Now put $m \in  {\mathbb N}$.
   The moduli spaces of determinants of packing lattices of  Minkowski balls and  Minkowski domains have the form
\begin{equation}
 \Delta(p,\sigma)_{2^m D_p} = 4^m (\tau + \sigma)(1 + \tau^{p})^{-\frac{1}{p}}
  (1 + \sigma^p)^{-\frac{1}{p}}.
   \end{equation}
 \end{proposition}
 

\subsection{Moduli spaces for covering by Minkowski balls
}
Every admissible lattice of $D_p$ containing 3 pairs of points on the boundary of $D_p$ defines a hexagon inscribed in $D_p$.
We denote such a hexagon as $ {\mathcal{H}} _p $, and call it a hexagon of the admissible lattice, briefly: {\it al-hexagon}.

\begin{proposition}
 \label{msa}
The set of areas of the entire family of al-hexagons  $ {\mathcal{H}}_p, (1 <  p < \infty)$ inscribed in Minkowski balls $D_p$ is parameterized by the function $A(\sigma, p)$ over the domain (\ref{msd}) and defining the corresponding moduli space
   $\mathbf{A}: $
 \begin{equation}
 \label{Ap}
 A(\sigma, p) = 3(\tau + \sigma)(1 + \tau^{p})^{-\frac{1}{p}}
  (1 + \sigma^p)^{-\frac{1}{p}}, 
 \end{equation}
\end{proposition}

\subsection{Moduli space of areas of al-hexagons inscribed in unit Minkowski balls}

From Proposition \ref{msa} we have:
 
\begin{corollary}
\label{msh}
 The moduli space of areas of al-hexagons inscribed in unit Minkowski balls has the form (\ref{Ap}).
\end{corollary}

\subsection{Moduli spaces of areas of al-hexagons circumscribed to unit two-dimensional Minkowski balls}

\begin{proposition}
 \label{msch}
The set of areas of the entire family of ( $ {\mathcal{H}}_p, (1 <  p < \infty)$ circumscribed to unit Minkowski balls has the form:
 \begin{equation}
{\mathbf {AC}}(p,\sigma)_{2 D_p} = 4 (\tau + \sigma)(1 + \tau^{p})^{-\frac{1}{p}}
  (1 + \sigma^p)^{-\frac{1}{p}}, 
 \end{equation}
in the same domain (\ref{msd}).
\end{proposition}

\section{Extreme functions and solution of optimization problems}

The existence of these extremal functions is based on the fact that the infimum and supremum in
 formulas    (\ref{cd}), (\ref{ham}), (\ref{ha})   are achieved.

\subsection{Extreme function of Minkowski-Davis problem}

 For  Minkowski-Cohn moduli space (\ref{mcms}) of admissible lattices of Minkowski balls we have

Let
\[
 {\Delta^{(0)}_p} = \Delta(p, {\sigma_p}) =  \frac{1}{2}{\sigma}_{p},   \; {\sigma}_{p} = (2^p - 1)^{1/p},
\]
\[
 {\Delta^{(1)}_p}  = \Delta(p,1) = 4^{-\frac{1}{p}}\frac{1 +\tau_p }{1 - \tau_p},   \;  2(1 - \tau_p)^p = 1 + \tau_p^p,  \;  0 \le \tau_p < 1.  
\]
The extremal function of critical lattices of Minkowski balls  has the form (\cite{GGM:PM}):

 \begin{equation}
 \label{efmb}
\Delta(D_p) = \left\{
                   \begin{array}{lc}
    \Delta(p,1), \; 1 < p \le 2, \; p \ge p_{0},\\
    \Delta(p,\sigma_p), \;  2 \le p \le p_{0};\\
                     \end{array}
                       \right.
                          \end{equation}
here $p_{0}$ is a real number that is defined unique by conditions
$\Delta(p_{0},\sigma_p) = \Delta(p_{0},1),  \;
2,57 < p_{0}  < 2,58, \; p_0  \approx 2.5725 $

\subsection{Extremal functions for packing  of Minkowski balls and domains}

\begin{theorem}Let $m \in {\mathbb N}$ (see \cite{gl1} for the case $m=1$).
The critical determinants of $2^m$ {\it  dubling}
  balls $D_p$ have a representation of the form
\begin{equation}
\label{cdde0}
 {\Delta^{(0)}_p}(2^m D_p) = \Delta(p, {\sigma_p})_{2^m D_p} =  2^{2 m-1}\cdot {\sigma}_{p},   \; {\sigma}_{p} = (2^p - 1)^{1/p},
 \end{equation}
 \begin{equation}
\label{cdde1}
 {\Delta^{(1)}_p}(2^m D_p)  = \Delta(p,1)_{2^m D_p} = 4^{m - \frac{1}{p}}\frac{1 +\tau_p }{1 - \tau_p},   \;  2(1 - \tau_p)^p = 1 + \tau_p^p,  \;  0 \le \tau_p < 1. 
 \end{equation}
  And these are the determinants of the sublattices of index $2^m$ of the critical lattices of the corresponding  balls $D_p$.
  
  The extremal functions for packing  of Minkowski balls and domains have the forms:
 \begin{equation}
 \label{efmbd}
\Delta(2^m D_p) = \left\{
                   \begin{array}{lc}
    \Delta^{(1)}_p (2^m D_p), \; 1 <  p \le 2, \; p \ge p_{0},\\
    \Delta^{(0)}_p (2^m D_p), \;  2 \le p \le p_{0};\\
                     \end{array}
                       \right.
                       \end{equation}
\end{theorem}

\subsection{Extremam function of covering by  Minkowski balls}

In this case the extremal function is unknown. It is possible to construct its approximations \cite{gl}.

\begin{problem*}
\label{mbf}
What are extremal function and covering constants for covering problem?
\end{problem*}

\subsection{Extremal functions for  hexagons of  minimal areas inscribed in and circumscribed to unit two-dimensional Minkowski balls}

From Proposition \ref{msa}, Corollary \ref{msh} and Proposition \ref{msch}    for  hexagons of  minimal areas inscribed in  two-dimensional Minkowski balls, we have thefollowing extremal function:

\begin{theorem}
\label{ma}
\begin{equation}
\label{ema}
\min({ A}({\mathcal {IH}}_p)(D_p)) = \left\{
                   \begin{array}{lc}
    3\cdot4^{-1/p} \frac{1 +\tau_p }{1 - \tau_p}, \; 1 < p \le 2, \; p \ge p_{0},\\
    \frac{3}{2}{\sigma}_{p}, \;  2 \le p \le p_{0};\\
                     \end{array}
                       \right.
                       \end{equation}
\end{theorem}

\begin{remark}
   From Propositions 2 and 5 it follows that the moduli spaces
 of packing lattices and hexagons circumscribed around 
Minkowski balls coincide. Consequently, their extremal
 functions also coincide. This function is a specialization of function (\ref{efmbd}) for $m = 1$.
 For completeness, below we present its form.
\end{remark}

\begin{theorem}
\label{ch}
  Below $ 0 \le \tau_p < 1,$ \\
\begin{equation}
\label{ech}
   \min({ A}({\mathcal {CH}}_p)(D_p)) = \left\{
                   \begin{array}{lc}
    4^{1 - \frac{1}{p}}\frac{1 +\tau_p }{1 - \tau_p},   \;  2(1 - \tau_p)^p = 1 + \tau_p^p,  \; 1 < p \le 2, \; p \ge p_{0},\\
     2 \cdot {\sigma}_{p},   \; {\sigma}_{p} = (2^p - 1)^{1/p},, \;  2 \le p \le p_{0};\\
                     \end{array}
                       \right.
   \end{equation}
\end{theorem}

{\bf Acknowledgements.}{ This work was carried out at Institute of Mathematics and Informatics (IMI) of the BAS.
 I thank  the Bulgarian Academy of Sciences, the IMI BAS for support.  The author is deeply grateful to P. Boyvalenkov and  V. Drensky for their support.}\\



 

\end{document}